\documentstyle[12pt,geompsfi]{article}
\font\tenmsb=msbm10 scaled \magstep1
\font\sevenmsb=msbm8
\font\fivemsb=msbm6
\newfam\msbfam
\textfont\msbfam=\tenmsb
\scriptfont\msbfam=\sevenmsb
\scriptscriptfont\msbfam=\fivemsb
\def\Bbb{\fam\msbfam}

\def\Z{{\Bbb Z}}

\def\R{{\Bbb R}}
\def\C{{\Bbb C}}

\newtheorem{thm}{Theorem}[section]

\newtheorem{prop}[thm]{Proposition}

\newtheorem{lemma}[thm]{Lemma}
\newtheorem{cor}[thm]{Corollary}

\newcommand{\rs}{\mbox{$ \widehat{ \C } $}}

\newcommand{\vac}{\mbox{\O}}
\newcommand{\qed}{\nopagebreak \begin{flushright}
	\rule{2mm}{2.5mm} \end{flushright}}

\newcommand{\bdry}{\partial}
\newcommand{\cl}{\overline}
\newcommand{\implies}{\mbox{$\Longrightarrow$}}

\newcommand{\gap}{\vspace{5pt}}
\newcommand{\pf}{\noindent {\bf Proof: }}

\newcommand{\rmk}{\gap \noindent {\bf Remark: }}

\newcommand{\Ins}{\mbox{Ins}}
\newcommand{\Out}{\mbox{Out}}
\newcommand{\Int}{\mbox{Int}}

\begin{document}

\title{Rational maps whose Fatou components are Jordan domains }

\author{Kevin M. Pilgrim \footnotemark}

\maketitle

\begin{abstract}
We prove: If $f(z)$ is a critically finite rational map which has
exactly two critical points and which is not conjugate to a
polynomial, then the boundary of every Fatou component of $f$ is a
Jordan curve.  If $f(z)$ is a hyperbolic critically finite rational
map all of whose postcritical points are periodic, then there exists a
cycle of Fatou components whose boundaries are Jordan curves.  We give
examples of critically finite hyperbolic rational maps $f$ with the
property that on the closure of a Fatou component $\Omega$ satisfying
$f(\Omega)=\Omega$, $f|_{\bdry \Omega}$ is not topologically conjugate
to the dynamics of any polynomial on its Julia set.
\end{abstract}

\footnotetext{Research supported by a National Need Fellowship and NSF
Grant DMS-9301502.  The author would like to thank J. Kahn and
C. McMullen for useful conversations.
Research at MSRI is supported in part by NSF grant no. DMS-9022140.
Address: c/o Mathematical Sciences
Research Institute, 1000 Centennial Drive, Berkeley, CA 94720, USA.
E-mail: pilgrim@msri.org}

\section{Introduction}

A {\it rational map} $f(z)=p(z)/q(z)$ where $p$ and $q$ are relatively
prime complex polynomials determines a holomorphic map of the Riemann
sphere $\rs $ to itself, and so defines a holomorphic dynamical
system.  The {\it Fatou set} $J(f)$ is the set of those $z \in \rs $
such that there exists a neighborhood $U$ of $z$ on which the iterates
$\{f^{n}|_{U}\}$ form a normal family of holomorphic functions.  The
complement $J(f)$ of the Fatou set is called the {\it Julia set}.  We
shall assume throughout that the degree $d$ of $f$ is larger than one.
The set of points $C(f)$ where the derivative of $f$ vanishes is the
set of {\it critical points} of $f$; these are the points where the
local degree of the map is larger than one.  Counted with
multiplicity, there are $2d-2$ critical points.

A characteristic feature of the dynamics of iterated rational maps is
that the behavior of the finite set of critical points under iteration
strongly influences the dynamics of the map on the entire sphere.  For
examples of this, and basic definitions, see e.g. \cite{MIL1}.  The
{\it postcritical set} $P(f)$ of a rational map $f$ is defined by
$P(f)=\cl{ \cup_{n>0} f^{n}(C(f)) }$.  The map $f$ is said to be {\it
critically finite} if $P(f)$ is finite.  The map $f$ is said to be
{\it hyperbolic} if $P(f) \cap J(f) = \vac$.

It is known that the Julia set of a critically finite map is
connected, and hence every Fatou component is an open disc.  The
boundaries of Fatou components for critically finite maps are known
also to be locally connected.  The boundaries of the Fatou components,
however, need not be Jordan curves.  A {\it Jordan domain} is a
component of the complement of a Jordan curve in $S^{2}$.

In this paper we prove

\begin{thm}
\label{thm:mainthm}
Let $f(z)$ be a critically finite rational map with exactly two critical
points, not counting with multiplicity.  Then exactly one of the
following possibilities holds:

\begin{itemize}
\item  $f$ is conjugate to $z^d$ and the Julia set of $f$ is a Jordan
curve, or

\item  $f$ is conjugate to a polynomial of the form $z^{d}+c, c \neq
0$, and the Fatou component corresponding to the basin of infinity
under a conjugacy is the unique Fatou component which is not a Jordan
domain, or

\item  $f$ is not conjugate to a polynomial, and every Fatou
component is a Jordan domain.

\end{itemize}
\end{thm}

Since a quadratic rational map has exactly two simple critical
points, the hypothesis of Theorem \ref{thm:mainthm} are satisfied for all
postcritically finite quadratic rational maps.  Theorem \ref{thm:mainthm}
confirms what had been experimentally observed in computer studies:
that for many critically finite quadratic rational maps which are not
polynomials, every Fatou component has Jordan curve boundary; see \cite{MIL2}.

\begin{thm}
\label{thm:isone}
Let $f$ be a hyperbolic critically finite rational map. If every
postcritical point of $f$ is periodic, then there exists at least one
cycle of Fatou components of $f$ consisting of Jordan domains, and
every Fatou component which maps onto an element of this cycle is also
a Jordan domain.
\end{thm}

\begin{cor}
\label{cor:samecycle}
If $P(f)$ consists of a single superattracting cycle, then
every Fatou component of $f$ is a Jordan domain.
\end{cor}

A proof of Theorem \ref{thm:mainthm} in the special case of certain
critically finite quadratic maps appeared in \cite{rees:parami},
Section 5.4, as an ingredient in the classification of quadratic
rational maps.  The argument given assumes the following
fact, specialized to the quadratic case:

\noindent{\bf Invariance Condition:} {\it Let $f$ be a critically
finite rational map with exactly two critical points.  Let $ \Omega $
be a periodic Fatou component of period $p$ for which
$f^{p}|_{\Omega}$ is conjugate to $z \mapsto z^{m}$.
Then for every $x \in \bdry \Omega$, $|(f^{p}|_{\bdry
\Omega})^{-1}(x)| \leq m$.}

This fact need not hold for maps with three or more critical points:
we give a degree three example in Section \ref{section:examples}.
These examples show (with the above notation) that $f^{p}|_{\bdry
\Omega}$ need not be topologically conjugate to the dynamics of any
polynomial on its Julia set, even in the hyperbolic case.  In
addition, a Jordan curve in $\bdry \Omega$ need not have a preimage
under $f^{p}$ which is contained in $\bdry \Omega$.  This shows that, if
one associates to $\Omega$ an invariant lamination $L$ in the sense of
Thurston (see \cite{thurston:dynamics}), then $L$ may fail to
satisfy the condition of {\it gap invariance}.  Since the writing of
this paper Tan Lei and the author have obtained a description of these
examples in terms of a new kind of surgery which will be the subject
of a future paper; see \cite{kmp:tan:blow}.

The process of {\it tuning} is a way of combining the dynamics of a
rational map $f$ with the dynamics of a polynomial $p$.  An open
question is to find conditions on a critically finite map $f$ and
a critically finite polynomial $g$ for the tuning to be
combinatorially equivalent to a rational map in the sense of Thurston.
If $f$ has a periodic Fatou component with non-Jordan curve boundary,
it is generally believed that there exists some $g$ for which the
tuning has a topological obstruction to being combinatorially
equivalent to a rational map.  The converse is known to be false; see
\cite{ahmadi:thesis}, Theorem 5.11.1.

In Section \ref{section:background} we state known facts from the
theory of iterated rational maps which we use in the proofs of the
theorems.  In Section \ref{section:curves} we develop the main
technique used in the proof: an analysis of how Jordan curves in the
Julia set behave under backwards iteration.  In Section
\ref{section:proofs} we prove Theorems \ref{thm:mainthm} and
\ref{thm:isone}.  In Section \ref{section:examples} we give examples
of non-polynomial maps which fail the Invariance Condition stated above.

\section{Background}
\label{section:background}

\subsection{An important consequence of Montel's theorem}

We will use the following proposition to control how a Jordan domain behaves
under backwards iteration of a rational map.

\begin{prop}[Montel's Theorem]
\label{inF}
Let $f$ be a rational map and $U \subset \rs$ be a connected open set
whose complement contains at least three points.  If $f^{-p}(U) \supset
U$ for some $p>0$, then $U$ is in the Fatou set.
\end{prop}

{\bf Proof:} For all $n \geq 0$, the images $(f^{p})^{n}|_{U}(U)$
omit at least three points, and hence form a normal family of holomorphic
functions by Montel's theorem.  The iterates of $f|_{U}$ then form a
normal family, and so $U$ is in the Fatou set.
\qed

\subsection{Postcritically finite rational maps}

In this subsection we collect needed facts about critically finite
maps.  In particular, these maps have important expanding
properties. For the definition of orbifold, the canonical orbifold
associated to a critically finite map, and the definition of the
associated canonical orbifold Poincar\'e or Euclidean metric, see
\cite{MIL1}, Lemma 14.5.  This metric is behaves very much
like the usual Poincar\'e or Euclidean metric on $\rs - P(f)$.  Let
$Q(f)$ denote the set of postcritical points which eventually land on
cycles containing critical points. Then $Q(f) \subset F(f)$.  The
canonical orbifold metric $\rho$ is supported on $\rs - Q(f) $ and
lifts under $f$ to a metric $\tilde{\rho}$ on $\rs - f^{-1}Q(f)$.
With respect to the metric $\tilde{\rho}$ on $\rs - f^{-1}Q(f) $ and
the metric $\rho$ on $\rs - Q(f) $, the inclusion $\rs - f^{-1}(Q(f))
\hookrightarrow \rs - Q(f)$ is a strict contraction.  We then have

\begin{prop}
\label{prop:cf maps expand}
Let $f$ be a critically finite rational map.  Then $f$ is uniformly
expanding with respect to the canonical orbifold metric $\rho$ on the
complement of any open neighborhood of $f^{-1}Q(f)$.  In particular,
$f$ is uniformly expanding on $J(f)$ with respect to $\rho$.
\end{prop}

{\bf Proof:} See \cite{MIL1}, Theorem 14.4.

\begin{prop} Let $f$ be a critically finite rational map.  Then
\begin{enumerate}
\item a Fatou component contains at most one critical point;
\item every Fatou component is eventually periodic;
\item the Julia set is connected;
\item the Julia set is the whole sphere iff there are no periodic
critical points.
\end{enumerate}
\end{prop}

{\bf Proof:} In a periodic Fatou component containing two or more
critical points, at least one must have an infinite forward orbit,
which cannot happen if the map is critically finite.  This implies
(1).  To prove (2), we may appeal to Sullivan's No Wandering Domains
Theorem \cite{sullivan:qcI}, but in our case one can use expansion of
the orbifold metric to give a direct argument.  The point is that a
wandering Fatou component must avoid a neighborhood of $Q(f)$, and so
each iterate of $f$ contributes a definite factor of expansion on such
a component.  Hence in a wandering sequence of components, the
diameters of the components must tend to infinity with respect to the
canonical orbifold metric, which is impossible.  This proves (2).  It
now follows that every periodic Fatou component is a covering of an
open disc branched over at most one point, hence every Fatou component
is a disc and so (3) is proved.  To prove (4), if there are no
periodic critical points, then $Q(f)$ is empty and so $f$ expands the
orbifold metric at every point of the Riemann sphere.  The Julia set
of $f$ is thus the entire sphere.  Conversely, a periodic critical
point is always in the Fatou set.
\qed

\begin{prop}
\label{boundary_lc}
Let $f$ be a critically finite rational map and $\Omega $ a period $p$
Fatou component.  Then $\partial \Omega $ is locally connected.
\end{prop}

{\bf Proof:} The proof is virtually identical to the proof of the
well-known corresponding fact for subhyperbolic polynomials with
connected Julia set; see \cite{MIL1}, Theorem 17.5. The only
difference is that one uses the first return map $f^{p}$ restricted to
$\Omega $ in place of the polynomial restricted to its basin of
infinity.
\qed

\subsection{Riemann mappings and local connectivity}

A set $K \subset \C $ is said to be {\it full} if it is compact,
connected, and if its complement is nonempty and connected.  A full
set is said to be {\it nondegenerate} if it is not a point.

We will need the following

\begin{prop}
\label{prop:isjcurve}
Let $K$ be a full nondegenerate subset of $\C$ whose boundary is
locally connected.  Let $V$ be a bounded component of $\rs -
\partial K $.  Then

\begin{enumerate}
\item $V$ is a Jordan domain,
\item $\cl{V}$ and $\rs - V$ are closed discs, and
\item a Jordan curve in $K$ is contained in the closure of a unique
bounded component $U$ of $\rs - \partial K$.
\end{enumerate}
\end{prop}

The first conclusion is essentially the content of
\cite{DH2}, Section 2.4.3, where it is stated without proof.
For completeness, we give a proof.  In what follows, $\Delta$ denotes
the open unit disc $\{z: |z| < 1 \}$ and $\Sigma = \rs - \cl{\Delta}$.

\begin{thm}[Carath\'eodory]
Let $K$ be a full nondegenerate set in $ \C $.  Let $\phi : (\Delta, 0
) \rightarrow (\rs - K, \infty )$ be a Riemann map
uniformizing the complement of $K$ in \rs .  Then $\phi $ extends to a
continuous map $\overline{\phi}: \overline{\Delta} \rightarrow \rs $
if and only if $\partial K$ is locally connected, or if and only if
$K$ is locally connected.
\end{thm}

See \cite{MIL1}, Theorem 16.6 for the proof.

Let $\phi: (\Delta, 0 ) \rightarrow (U, z)$ be a Riemann map
uniformizing an open disc $U$.  For $t \in \R / \Z $ the {\it ray of
angle $t$ for $\phi $ }is the set $ \{ \phi( re^{2\pi i t} ) | r \in
[0, 1) \} $, and is denoted by $R_{t}$.  If $R_{t}$ has a unique limit
point $x$ in $\partial U$, the ray $R_{t}$ is said to {\it land} at
$x$.

\begin{thm}
\label{rays_separate}
Let $K$ be a full nondegenerate locally connected set in $\C $ and $U=
\rs - K$ .  Let $\phi : (\Delta, 0) \rightarrow (U, \infty) $ be a
Riemann mapping.  Suppose two distinct rays $R_{t}$ and $R_{t'}$ of
$\phi $ land at a common point $x$ of $\partial U$.  Then $x$ separates
$\partial K$ so that each component of the complement of the Jordan
curve $C=R_{t} \cup R_{t'} \cup \{x\}$ contains a nonempty component
of $\partial K - \{x\}$.

\end{thm}

{\bf Proof:} Since $\partial K$ is locally connected, $\phi$ extends
to a map $\overline{\phi }$ of the closed disc, by Carath\'eodory's
Theorem.  Suppose $C$ failed to separate $\partial K$ so that some
component of its complement did not contain points of $\partial K$.
Since $\phi$ is a homeomorphism and since rays cannot cross in $U$, by
relabelling $t$ and $t'$ we may assume that for the set $W=\{re^{2\pi
i s} | |r| < 1, s \in (t,t') \}, \cl{ \phi(W)} \cap K = \{x\} $.  So
$\cl{\phi}$ collapses the circular arc $(t,t') \subset S^{1}$ to the
point $x$.  But this contradicts the Theorem of F. and M. Riesz
\cite{caratheodory}, Volume II, Section 313.  \qed

{\bf Proof of Proposition \ref{prop:isjcurve}:} Since $K$ is full, its
boundary is connected.  Hence every bounded component of the
complement of $\partial K$ is an open disc.

Let $V$ be a bounded component of $\C - \partial K$.  Note
first that $\overline{V} \subset K $, since $K$ is full.   We first
show that $\partial V$ is locally connected.  Let $\phi : (\Sigma,
\infty) \rightarrow (\rs - K, \infty) $ be a Riemann map to
the complement of $K$ in \rs .  By Carath\'eodory's Theorem, $\phi $
extends to a continuous map on $\overline{\Sigma }$.  Consider the
equivalence relation on $S^{1}$ defined as follows: $ x \sim y$ iff
$\phi (x) = \phi (y) $.  For each equivalence class in $S^{1}$, form
its Euclidean convex hull in $\overline{\Delta }$.  The convex hulls
of any two distinct equivalence classes are disjoint.  Let $L$ be the
union of the convex hulls of equivalence classes.  Then $\phi $
extends to a map $\overline{\phi }: \cl{\Sigma} \cup L \rightarrow
\cl{\rs - K }$ by mapping the convex hull of any equivalence class
$[x]$ to $\phi(x)$.  Then $\partial K =
\overline{\phi }(S^{1}\cup  L)$ and $\overline{\phi }$ is a
homeomorphism from $\Sigma$ to $\rs -  K $.

Given $V$, let $x$ and $y$ be distinct points on $\partial V$.  Then
$\cl{\phi}^{-1}(x)$ and $\cl{\phi}^{-1}(y)$ are not separated in
$\cl{\Delta}$ by the preimage of any other point $z$ in $\partial K$
distinct from $x$ and $y$.  For otherwise there are rays $R_{t},
R_{t'}$ landing at $z$ such that $C=R_{t} \cup \{z\} \cup R_{t'}$ is a
simple closed curve in $\cl{ \rs - K}$ separating $x$ and $y$,
contradicting the fact that $x$ and $y$ lie in the boundary of a
single component of $\rs - \partial K$.

Let $X=\cl{\phi}^{-1}(\partial V)$, and let $\cal C \cal H (X)$ be the
Euclidean convex hull of $X$ in $\cl{\Delta}$.  Then $\cal C \cal H
(X) \supset X$, and since it is the convex hull of $X$ its boundary is
contained in $\partial X$. Its boundary is locally connected since the
boundary of any bounded convex set is locally connected; see
\cite{newman:plane_sets}, ch. 6 section 4.  Since $X$ is closed,
$\partial X \subset X$, and so $\cl{\phi}(\partial \cal C \cal H (X) )
\subset \partial V$.  Moreover, $\cl{\phi}(\partial \cal C \cal H (X)
) = \partial V$.  So $\partial V$ is the continuous image of a compact
locally connected set, and hence is locally connected, by Lemma 16.5
of \cite{MIL1}.  In particular, every point on $\partial V$ is the
unique limit point of some ray of $\psi $, where $\psi: \Delta
\rightarrow V$ is a Riemann map uniformizing $V$.

We now show that $\partial V$ is a Jordan curve.  If $\partial V$ is
not a Jordan curve, there is a point $p \in \partial V$ which is the
landing point of at least two rays for $\psi $.  The union of these
rays, together with the common landing point, gives a Jordan curve $C
\subset \overline{V}$ which separates $\partial V$ into at least two
pieces, one of which lies in the bounded component $W$ of $\rs - C$ by
Theorem \ref{rays_separate}. Then $\partial K
\cap W \neq \vac $, hence there exist points of $\rs - K$ in
$W$.  Since $\partial W \subset \overline{V} \subset K$, this implies
that $\partial W$ separates points in $\rs - K$, and hence
that $\rs - K$ is not connected, contradicting $K$ full.

Hence $V$ is an open disc with Jordan curve boundary.  The Schoenflies
theorem implies that $\cl{V}$ and $\rs - V$ are both homeomorphic
to closed discs.  If $C$ is any Jordan curve contained in $K$, let $W$
be the bounded component of $\rs - C$.  Then $W$ is contained
in the interior of $K$.  If $\cl{W}$ is not contained in the closure
of a unique bounded component $V$ of $\rs - \partial K$, then $W$
must contain points of $\partial K$, which is impossible.
\qed

\subsection{Topological propositions}

We will make extensive use of the following

\begin{prop}
\label{prop:at most one critical value}
If $f$ is a rational map and $U$ is a Jordan domain whose closure
contains at most one critical value $v$ of $f$, then every component $V$
of $f^{-1}(U)$ is also a Jordan domain.  If $v \in \bdry U$, then
$f|_{\cl{V}}: \cl{V} \rightarrow \cl{U}$ is a homeomorphism.
\end{prop}

\pf Suppose first that $v \in U$.  Let $U'=\cl{U}-\{v\}$.  Let $V'$ be
the unique component of $f^{-1}(U')$ which contains a point of $V$.
Then $f|_{V'}: V' \rightarrow U'$ is an unbranched covering of a
closed punctured disc, and so $V'$ is a closed punctured disc.  Hence
$\cl{V}$ is a closed disc, and so $V$ is a Jordan domain.

If $v \in \bdry U$, it suffices to show that $f|_{\cl{V}}$ is
injective.  Suppose otherwise.  Then there exist $x_{1}, x_{2} \in
\bdry V$ such that $f(x_{1}) = f(x_{2})=y \in \bdry U$.  Choose $v \in
V$, and let $\alpha_{i}, i=1,2$ be closed, embedded arcs in $\cl{V}$
whose interiors are disjoint open arcs in $V$ and which join $v$ to
$x_{i}$.  Then $f(\alpha_{1} \cup \alpha_{2})$ is a Jordan curve in
$\cl{U}$ which intersects $\bdry U$ in exactly one point.  Hence there
is a component $W$ of $V-(\alpha_{1} \cup \alpha_{2})$ such that
$f(\bdry W \cap \bdry V)=y$.  The set $\bdry W \cap \bdry V$ is not discrete,
since $x_{1} \neq x_{2}$.  But this is impossible since $f$ is
holomorphic on $\rs $.
\qed

\section{Jordan curves in $J(f)$ }
\label{section:curves}
\newcommand{\sign}{\mbox{sign}}

In this section we develop the techniques used in the proof of
Theorems \ref{thm:mainthm} and \ref{thm:isone}.
\gap

\noindent {\bf Convention:} Throughout this section, $f(z)$ will
denote either a postcritically finite rational map with at most one
critical point in $J(f)$, or an iterate of such a map.
\gap

Given an {\it oriented} Jordan curve $\gamma$ in the sphere, define
the {\it inside} of $\gamma$, denoted by $\Ins (\gamma)$, to be the
component of the complement of $\gamma$ lying to the left of $\gamma$,
and the {\it outside} $\Out (\gamma)$ of $\gamma$ to be the component
of the complement lying to the right of $\gamma $.

Let $\gamma \subset J(f)$ be an unoriented Jordan curve.  A {\it lift}
of $\gamma$ we define to be a Jordan curve $\eta
\subset f^{-1}(\gamma)$ such that $f|_{\eta}: \eta \rightarrow \gamma$
is a covering map.  If $\gamma \subset J(f)$ is oriented, a {\it lift}
of $\gamma$ is an unoriented lift $\eta$ of $\gamma$, equipped with an
orientation so that $f|_{\eta}: \eta \rightarrow \gamma$ is
orientation-preserving.  Since $f$ is an iterate of a map with at most
one critical point in $J(f)$, a component of $f^{-1}(\gamma)$ is
homeomorphic either to a Jordan curve, or to a one-point union of
Jordan curves.  Thus if $\deg (f)=d$, there are exactly $d$ lifts of
any Jordan curve in $J(f)$, counted with multiplicity equal to the
absolute value of the degree of the map $f|_{\eta}: \eta \rightarrow
\gamma$.  We denote the set of lifts of a Jordan curve $\gamma$ by
$f^{*}(\gamma)$.

Let $\Omega$ be a Fatou component of $f$.  It is convenient to replace
the map $f$ by a conjugate so that some point of $\Omega$ is the point at
infinity.  Then $\rs - \Omega$ is a full locally connected subset of
the plane which has empty interior if and only if $\Omega$ is the
unique Fatou component of $f$, which in turn holds if and only if $f(z)$ is
conjugate to a polynomial with a unique Fatou component which is the
basin of infinity.  We shall therefore make the additional assumption
that $f(z)$ is not conjugate to such a polynomial.

An oriented Jordan curve $\gamma$ is said to be {\it positively
(negatively) oriented with respect to $\Omega$} if $\Omega \subset
\Out (\gamma)$ (respectively $\Omega \subset \Ins (\gamma))$.  We also
say that the {\it sign} $\sign (\gamma)$ of $\gamma$ is positive
(negative) if it is positively (negatively) oriented.  We denote
by $\gamma$ both an oriented and an unoriented Jordan curve; in
the following we will explicitly mention which is meant.

\noindent {\bf Notation:}

\begin{itemize}

\item Let $\Gamma$ denote the set of unoriented Jordan curves in the Julia
set of $f$.

\item Let $A$ denote the set of closures of components of $\rs -
\Omega$.  By Proposition \ref{prop:isjcurve}, each element $a \in A$ is a
closed disc.  The proposition also implies that there is a
well-defined function $p_{A}: \Gamma \rightarrow A$ which
assigns to every $\gamma \in \Gamma$ the unique element $a \in
A$ for which $\gamma \subset a$.

\item Let $\Gamma^{\pm}$ denote the set of {\it oriented} Jordan
curves in $J(f)$.  Then $\Gamma^{\pm}=\Gamma^{+} \cup \Gamma^{-}$,
where $\Gamma^{+}$ ($\Gamma^{-}$) is the set of curves which are
postively (respectively negatively) oriented with respect to $\Omega$.
The function $p_{A}$ extends naturally to $\Gamma^{\pm}$ by forgetting
the sign and then applying $p_{A}$; the composition we denote again by
$p_{A}$.

\item Let $\Sigma \Gamma ^{\pm}$ denote the set of infinite sequences
$\{\gamma_{n}\}_{n=0}^{\infty}$ satisfying: $\gamma_{n} \in
\Gamma^{\pm}$, $\gamma_{n+1} \in f^{*}(\gamma_{n})$, and
$\gamma_{n+1}$ is equipped with an orientation so that
$f|_{\gamma_{n+1}}: \gamma_{n+1} \rightarrow \gamma_{n}$ is
orientation-preserving.  If $\sign(\gamma_{n}) \neq \sign
(\gamma_{n+1})$ we say that the sequence
$\{\gamma_{n}\}_{n=0}^{\infty}$ has a {\it sign change} between
$\gamma_{n}$ and $\gamma_{n+1}$.

\end{itemize}

The idea for the proofs of our theorems is the following: consider
the subset $S(\Omega)$ of $\Sigma \Gamma ^{\pm}$ consisting of
sequences $\{\gamma_{n}\}$ for which $\gamma_{0} \in \bdry
\Omega$.  We use Proposition \ref{inF}, {\it Montel's theorem}, and
the fact that the preimages of disjoint sets are disjoint to deduce
relationships between sign changes in elements of $S(\Omega)$ and the
value of the function $p_{A}$ at the terms where sign changes occur.
In the special cases where the hypothesis of the theorems are
satisfied, this in turn will yield information about $\bdry \Omega$.

It turns out, however, that the set $S(\Omega)$ is too large to be used in
this manner, so we introduce a smaller space which captures the
features in which we are interested.  For any finite collection of
disjoint unoriented Jordan curves in $J(f)$, there is at least one
curve which is {\it outermost} in the following sense: it is not
separated from $\Omega$ by any other curve in the collection.  It is
not unique, in general.  We define $S_{out}(\Omega)$ to be the set of
sequences $\{\gamma_{n}\}_{n=0}^{\infty} \in S(\Omega)$ such that the
following holds: given any two consecutive terms $\gamma_{n+1},
\gamma_{n}$ regarded as unoriented curves, $\gamma_{n+1}$ is
outermost among the collection of lifts of $\gamma_{n}$.  We then say
that $\gamma_{n+1}$ is {\it outermost} among lifts of $\gamma_{n}$.

Let $\{\gamma_{n}\}_{n=0}^{\infty} \in S_{out}(\Omega)$.  Fix some $n
\geq 0$.  If $\sign (\gamma_{n})$ is positive, set $U = \Ins
(\gamma_{n})$ and $V=\Ins(\gamma_{n+1})$; otherwise let $U=\Out
(\gamma_{n})$ and $V = \Out (\gamma_{n+1})$.  In both cases, $U
\subset \rs - \Omega$.  Let $R$ be the the unique component of
$f^{-1}(U)$ for which $\gamma_{n+1} \subset \bdry R$.

\begin{figure}
\begin{center}~
\psfig{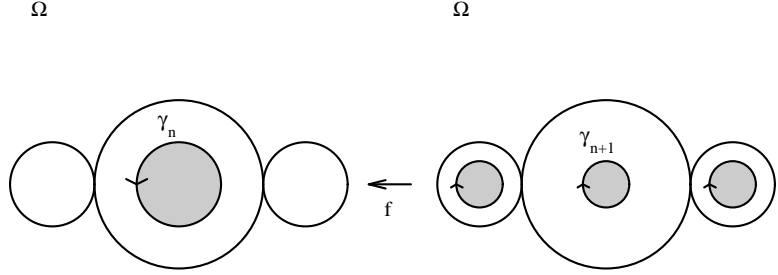}
\end{center}
\caption{$U=\Ins(\gamma_{n})$ is the shaded region on the left.  $R$
is the complement of the shaded regions on the right.  Only one
outermost lift $\gamma_{n+1} \subset \bdry R$ is labelled.
\label{fig:sign change} }
\end{figure}

\begin{figure}
\begin{center}~
\psfig{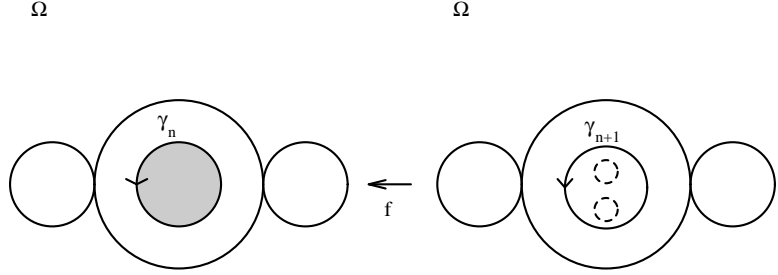}
\end{center}
\caption{$U=\Ins(\gamma_{n})$ is the shaded region on the left.  $R$
is region on the right bounded by $\gamma_{n+1}$ and the two dashed
Jordan curves, which are not outermost.
\label{fig:no sign change} }
\end{figure}

\begin{prop}[Sign changes]
\label{prop:sign_changes}
Given the hypotheses and notation in the preceding paragraph:
\begin{enumerate}

\item There is a sign change between $\gamma_{n}$ and
$\gamma_{n+1}$ if and only if $V \supset R \supset \Omega$ (Figure
\ref{fig:sign change}). Furthermore, if there is a sign change,
\[\bdry R =
\bigcup_{ \stackrel{ \delta \in f^{*}(\gamma_{n})}
                   { \delta \; \mbox{{\tiny outermost}}}
        } \delta . \]

\item There is no sign change between $\gamma_{n}$ and
$\gamma_{n+1}$ if and only if $R \subset V \subset \rs - \Omega$
(Figure \ref{fig:no sign change}).  Furthermore, if there is no sign change,
\[ f^{-1}(U) \subset \bigcup_{\stackrel{\delta \in
f^{*}(\gamma_{n})}{\delta \; \mbox{ {\tiny outermost}}}} V_{\delta}
\; \; \subset \rs - \Omega , \]
where $V_{\delta}= \Ins(\delta)$ if $\delta$ is positive and
$V_{\delta}=\Out (\delta)$ otherwise.

\item All outermost lifts of $\gamma_{n}$ have the same sign.
\end{enumerate}
\end{prop}

\pf  That $R \subset V$ is clear in both cases since $V$ is a disc and
$\gamma_{n+1}=\bdry V \subset \bdry R$.
\gap

\noindent {\bf 1:} Since $\bdry V=\gamma_{n+1}$ is outermost among
lifts of $\gamma_{n}$, there are no components of $\bdry R$ separating
$\bdry V$ from $\Omega$.  Since $\bdry R \subset \rs - \Omega$ and
there is a sign change, we must have $R \supset \Omega$.  Since every
other component of $f^{-1}(U)$ is disjoint from $R$, $R \supset
\Omega$, and $\bdry R \subset \rs - \Omega$, every other component of
$f^{-1}(U)$ is contained in $\rs - \Omega$.  Hence every other
component of $f^{-1}(U)$ is separated from $\Omega$ by a lift of
$\gamma_{n}$ which is contained in $\bdry R$.  Hence any outermost
lift of $\gamma_{n}$ must be contained in $\bdry R$. Since $R \supset
\Omega$, every Jordan curve in $\bdry R$ must be outermost among lifts
of $\gamma_{n}$.  The other implication is then clear.
\gap

\noindent {\bf 3:}  By Part 1, if there is a sign change between
$\gamma_{n}$ and $\gamma_{n+1}$, then all outermost lifts have the
same sign since they comprise $\bdry R \supset \Omega$.  Hence all
outermost lifts must have the same sign.
\gap

\noindent {\bf 2:} By Part 1, all outermost lifts of $\gamma_{n}$ have
the same sign.  Since $U \subset \rs - \Omega$ and there is no sign
change between $\bdry U = \gamma_{n}$ and $\bdry V = \gamma_{n+1}$, $V
\subset \rs - \Omega$.  Applying this to the collection of outermost
lifts and using the definition of outermost proves the second
assertion.
\qed

As a corollary to the previous proposition, we have

\begin{prop}[$\Omega$ fixed iff no sign changes]
\label{prop:fixed iff no sign changes}
The following are equivalent:
\begin{enumerate}
\item $f(\Omega)=\Omega$;

\item for every sequence $\{\gamma_{n}\}_{n=0}^{\infty} \in
S_{out}(\Omega)$, $\sign(\gamma_{1})=\sign(\gamma_{0})$.

\item for every sequence $\{\gamma_{n}\}_{n=0}^{\infty} \in
S_{out}(\Omega)$, there are no sign changes;

\end{enumerate}
\end{prop}

\pf That 3 $\implies 2$ is obvious.

\noindent {\bf 1 $\iff$ 2:}  By the previous proposition, $2$ holds if
and only if for every $a \in A$, $f^{-1}(a) \subset \rs - \Omega$.
Since $\rs - \Omega = \cl{\cup_{a \in A}a}$, we have that 2 holds if and
only if $f^{-1}(\rs - \Omega) \subset \rs - \Omega$, which in turn
holds if and only if $f(\Omega) \subset \Omega$.  Since $\Omega$ is a
Fatou component of $f$, this holds if and only if $f(\Omega)=\Omega$.

\noindent {\bf 1 $\implies$ 3:}  If there is a sign change between
$\gamma_{n}$ and $\gamma_{n+1}$, let $U, R$ be as in Part 1 of the
previous proposition.  Then $f^{-1}(U) \supset
\Omega$.  Since $U \cap \Omega = \vac$, we have $f^{-1}(U) \cap
f^{-1}(\Omega) = \vac$.  Hence $f^{-1}(\Omega) \subset \rs - f^{-1}(U)
\subset \rs - \Omega$.  But then $f^{-1}(\Omega) \not\supset \Omega$,
therefore $f(\Omega) \not = \Omega$.
\qed

The next proposition relates the dynamics of curves in the boundary of
a Fatou component $\Omega$ for which $f(\Omega)=\Omega$ to the
dynamics inside $\rs - \cl{\Omega}$.

\begin{prop}[$\Omega$ fixed]
\label{prop:Omega fixed}
Suppose $\Omega$ is forward-invariant under $f$.  Let $E \subset
\Int (a)$ be a nonempty subset, and suppose $f^{-1}E \subset
\cup_{i=1}^{i=r}b_{i}$, where  the $b_{i}$'s are distinct elements of
$A$ and $b_{i}\cap f^{-1}(E) \neq \vac$ for each $i$. Let $\gamma_{0}$ be the
positively oriented boundary of $a$.  If $\gamma_{1} \in
f^{*}(\gamma_{0})$, then $p_{A}(\gamma_{1})=b_{i}$ for some $b_{i}$.
\end{prop}

\rmk If $f(\Omega)=\Omega$, a Jordan curve $\gamma \subset \bdry
\Omega$ need not have a lift which is contained in $\bdry \Omega$; see
the examples in Section \ref{section:examples}.

\pf We may assume that $\gamma_{1}$ is outermost among lifts of
$\gamma_{0}$.  By Part 2 of Proposition \ref{prop:sign_changes}, {\it Sign
changes}, $f^{-1}(E) \subset f^{-1}(\Ins(\gamma_{0}))
\subset \rs - \Omega$.  Hence for $b \in A$, $f^{-1}(E) \cap b \neq \vac \iff
f^{-1}(\Ins( \gamma_{0})) \cap b \neq \vac$.  Since
$f^{-1}(\Ins(\gamma_{0})) \subset \rs - \Omega$,
$f^{-1}(\Ins(\gamma_{0})) \cap b \neq \vac $ if and only if there is a
collection $\gamma_{1}^{1}, \gamma_{1}^{2}, ..., \gamma_{1}^{k}$ of
outermost lifts of $\gamma_{0}$ such that $f^{-1}(\Ins(\gamma)) \cap b
\subset \cup_{j=1}^{k}\Ins(\gamma_{1}^{j})$.  This proves the proposition.
\qed

The next proposition refines the conclusion of the previous one in the
case where the topology of the map $f$ is simple.  Note that the
hypothesis is on the preimage of $\Out(\gamma_{0})$, which contains $\Omega$.

\begin{prop}[$\Omega$ fixed plus disc preimage]
\label{prop:Omega fixed plus disc preimage}
Given the hypothesis in Proposition \ref{prop:Omega fixed}, $\Omega$
{\it fixed}, suppose further that the component $V$ of the preimage of
$\Out (\gamma_{0})$ containing $\Omega$ is a Jordan domain.  Then
$f^{-1}E \subset b$ for a unique $b \in A$, and $\Omega$ is a Jordan
domain if and only if $a=b$.
\end{prop}

\pf The boundary of $V$ is the unique outermost preimage of
$\gamma_{0}$, since $V$ is a Jordan domain.  Thus $f^{-1}E$ is
contained in a unique $b \in A$, by Proposition \ref{prop:Omega
fixed}, {\it $\Omega$ fixed}.  If $\Omega$ is already a Jordan domain
the statement is trivially satisfied; the other direction follows from
the fact that if $a=b$, then $\bdry V \subset a$ and so $V \supset
\Out (\gamma_{0})=f(V)$.  But then $V$ is in the Fatou set, by
Proposition \ref{inF}, {\it Montel's Theorem}. Since $\bdry V$ is a
Jordan curve in $J(f)$, $V=\Omega$ is a Jordan domain.
\qed

The next sequence of propositions treat the case when there are sign
changes, i.e. when $f(\Omega) \neq \Omega$.

\begin{prop}[$\Omega$ not fixed]
\label{prop:Omega not fixed}  Suppose $f(\Omega) \subset a_{\Omega}
\in A$.

\begin{enumerate}

\item Let $\eta_{0} = \bdry a_{\Omega}$, equipped with positive
orientation.  Then there exists a negatively oriented lift
$\eta_{1}$ of $\eta_{0}$ which is outermost among lifts of
$\eta_{0}$.

\item If $\{\gamma_{n}\}_{n=0}^{\infty}$ has a sign change between
$\gamma_{n}$ and $\gamma_{n+1}$, then

\begin{enumerate}

\item $p_{A}(\gamma_{n})= a_{\Omega}$,

\item there is a lift $\gamma_{n+1}'$ of $\gamma_{n}$
which is outermost among lifts of $\gamma_{n}'$ such that
$p_{A}(\gamma_{n+1}')= a_{\Omega}$.

\end{enumerate}
\end{enumerate}
\end{prop}

\pf
\noindent {\bf 1:} Since $f(\Omega) \subset \Ins(\eta_{0})=a_{\Omega}$,
$f^{-1}(\Ins(\eta_{0})) \supset \Omega$.  Hence by Part 1 of Proposition
\ref{prop:sign_changes}, {\it Sign changes}, there is a negatively
oriented lift $\eta_{1}$ of $\eta_{0}$ which is outermost among lifts
of $\eta_{0}$. This proves the first assertion.
\gap

\noindent {\bf 2(a):}We argue by contradiction. We may assume that
$\gamma_{n}$ is postively oriented (otherwise, replace $\Ins$ with
$\Out$ in what follows).  If $\gamma_{n}$ is contained in
some $a \neq a_{\Omega}$, then $\Ins(\gamma_{n}) \cap \Ins(\eta_{0})
= \vac$, hence $f^{-1}(\Ins(\gamma_{n})) \cap f^{-1}(\Ins(\eta_{0}))
= \vac$.  Since $\Omega \subset f^{-1}(\Ins(\eta_{0}))$,
$f^{-1}(\Ins(\gamma_{n}))$ cannot contain $\Omega$.  By Proposition
\ref{prop:sign_changes}, {\it Sign changes}, this implies that
$\sign(\gamma_{n+1})=\sign(\gamma_{n})$.
\gap

\noindent {\bf 2(b):} We now prove the remaining assertion by
contradiction.  Again, we may assume that $\gamma_{n}$ is positively
oriented.  Suppose no outermost lift of $\gamma_{n}$ is contained in
$a_{\Omega}=p_{A}(\gamma_{n})$.  Let $R$ be the component of
$f^{-1}(\Ins(\gamma_{n}))$ which contains $\Omega$.  Then by
Proposition \ref{prop:sign_changes}, {\it Sign changes}, $\bdry R$
forms the collection of outermost lifts of $\gamma_{n}$.  If no such
lift is contained in $p_{A}(\gamma_{n})=a_{\Omega}$, then $R \supset
a_{\Omega} \supset \Ins (\gamma_{n})=f(R)$.  Hence by Proposition
\ref{inF}, {\it Montel's theorem}, $R$ is contained in the Fatou set
of $f$, and this is impossible since $\bdry a_{\Omega} \subset R
\subset F(f)$ while at the same time $\bdry a_{\Omega} \subset J(f)$.
\qed

As before, we now refine the conclusion of the previous proposition
in the case when the topology of the map is simple.

\begin{prop}[$\Omega$ not fixed plus disc preimages]
\label{prop:Omega not fixed plus disc preimages}
Suppose for all oriented Jordan curves $\gamma \in J(f)$, equipped
with positive orientation relative to $\Omega$, every component of
$\Ins (\gamma)$ is a Jordan domain.  Then
\begin{enumerate}

\item If $\sign (\gamma_{n}) \neq \sign (\gamma_{n+1})$, then
$p_{A}(\gamma_{n})=p_{A}(\gamma_{n+1})= a_{\Omega}$, i.e. any sign
changes are concentrated in $a_{\Omega}$.

\item For any sequence $\{\gamma_{n}\}_{n=0}^{\infty} \in S_{out}(\Omega)$,
if $p_{A}(\gamma_{n}) \neq a_{\Omega}$, there are no sign
changes after the $n$th term;

\item $(\cup_{n>0}f^{n}\Omega)-\Omega \subset a_{\Omega}$;

\item If $\eta_{0}$ denotes the positively oriented boundary of
$a_{\Omega}$, then there exists a unique outermost negatively oriented lift
$\eta_{1}$ of $\eta_{0}$ which is contained in $a_{\Omega}$.

\item $f^{-1}\Omega \subset a_{\Omega}$.

\end{enumerate}
\end{prop}

\pf

\noindent {\bf 1:} By Part 1 of Proposition \ref{prop:sign_changes}, {\it Sign
changes}, there exists a component $R$ of $f^{-1}(\Ins (\gamma_{n}))$
which contains $\Omega$.  $R$ is a Jordan domain, by the hypothesis.
Hence there is a unique outermost lift $\gamma_{n+1}$.  The first
statement then follows from Part 2 of Proposition \ref{prop:Omega not fixed},
{\it $\Omega$ not fixed}.
\gap

\noindent {\bf 2:} Let $\{\gamma_{n}\}_{n=0}^{\infty} \in
S_{out}(\Omega)$ be any sequence containing a sign change and suppose
$p_{A}(\gamma_{n}) \neq a_{\Omega}$.  Let $k$ be the smallest postive
integer such that there is a sign change between the $(n+k)$th and the
$(n+k+1)$st term.  We may assume that $\gamma_{n+k}$ is positively
oriented and $\gamma_{n+k+1}$ negatively oriented (otherwise, replace
$\Ins$ with $\Out$ in what follows).  By hypothesis, and
the assumption that there are no sign changes until the $(n+k+1)$st
term, for $i=1,2,...,k-1$, the unique component $R_{n+i+1}$ of
$f^{-i}(\Ins(\gamma_{n}))$ containing $\gamma_{n+i+1}$ is a Jordan
domain contained in $\rs - \Omega$ whose boundary is $\gamma_{n+i+1}$.
Since there is a sign change between $\gamma_{n+k}$ and
$\gamma_{n+k+1}$, by Part 1 above, $p_{A}(\gamma_{n+k+1}) = a_{\Omega}$.
Let $R_{n+k+1}$ be the unique component of $f^{-1}(R_{n+k})$ whose
boundary is $\gamma_{n+k+1}$.  Then since there is a sign change, by
Part 1 of Proposition \ref{prop:sign_changes}, {\it Sign changes},
$R_{n+k+1} \supset \Omega$.  Hence $R_{n+k+1} \supset
\Ins(\gamma_{n})=f^{k+1}(R_{k+1})$.  Proposition \ref{inF}, {\it Montel's
theorem}, then implies that $R_{n+1}$ is in the Fatou set, which is
impossible.
\gap

\noindent {\bf 3:}  A consequence of Part 2 is the following: if
$\gamma_{n}$ is a positively oriented element of the sequence
$\{\gamma_{n}\}_{n=0}^{\infty} \in S_{out}(\Omega)$, and if
$p_{A}(\gamma_{n}) = a \neq a_{0}$, then every lift of $\gamma_{n}$ is
outermost.  For since there are no sign changes, $f^{-1}(\Ins(\gamma_{n}))
\cap \Omega = \vac$, by Part 2 of Proposition \ref{prop:sign_changes},
{\it Sign changes}.  Every component of $f^{-1}(\Ins(\gamma_{n}))$ is a Jordan
domain in $\rs - \Omega$, by hypothesis.  Hence no boundary component of
$f^{-1}(\Ins(\gamma_{n}))$ separates another boundary component from
$\Omega$, and hence every lift of $\gamma_{n}$ is outermost.

This observation has the following strong consequence: if $a \neq
a_{\Omega}$, and if $\gamma_{0}$ is the positively oriented boundary
of $a$, then $\cup_{n \geq 0} f^{-n}(\Ins(\gamma_{0})) \cap \Omega =
\vac$.  We prove this by contradiction.  Let $n \geq 0$ be the
smallest postive integer for which $f^{-(n+1)}(a) \cap \Omega \neq
\vac$.  By induction and the result of the preceding paragraph, for
every $ 1 \leq i \leq n$, every component of
$f^{-i}(\Ins(\gamma_{0}))$ is a Jordan domain in $\rs - \Omega$ whose
boundary $\gamma_{i}$ is outermost among lifts of $f(\gamma_{i})$.  If
$f^{-(n+1)}(\Ins(\gamma_{0})) \cap \Omega
\neq \vac$, there is a finite sequence $\{\gamma_{i}\}_{i=0}^{n+1}$
of curves in $J(f)$ such that $\gamma_{i+1} \in f^{*}(\gamma_{i})$,
$\gamma_{i+1}$ is outermost among lifts of $\gamma_{i}$ for all $i
\leq n$, and $\sign(\gamma_{n+1}) \neq \sign(\gamma_{n})$ (by Part 1
of Proposition \ref{prop:sign_changes}), {\it Sign changes}.  But this
violates the conclusion of Part 2.

The result in the preceding paragraph implies that no component in the
forward orbit of $\Omega$ can intersect an element of $A$ which is
not $a_{\Omega}$.
\gap

\noindent {\bf 4:} By Part 1 of Proposition \ref{prop:Omega not fixed}, {\it
$\Omega$ not fixed}, existence is clear.  Uniqueness follows since the
unique component of the preimage of $\Ins(\eta_{0})$ containing
$\Omega$ is a Jordan domain.

\noindent {\bf 5:} By the previous step, $f^{-1}(\Ins(\eta_{0}))$ is a
Jordan domain containing $\Omega$ whose boundary is contained in
$a_{\Omega}$.  Hence $f^{-1}(\Omega) \subset f^{-1}(\Out(\eta_{0}))
\subset \rs - f^{-1}(\cl{\Ins(\eta_{0})}) \subset a_{\Omega}$.
\qed

\section{Proofs of the theorems}
\label{section:proofs}

\subsection{Proof of Theorem 2}
\label{perthm}

If $f$ is hyperbolic, there are no critical or postcritical points in
$J(f)$, so we may apply the analysis in Section \ref{section:curves}.
The proof of Theorem \ref{thm:isone} is then essentially a
straightforward application of Proposition \ref{prop:Omega fixed plus
disc preimage}, {\it $\Omega$ fixed}.

Choose arbitrarily an element $x \in P(f)$.  There is a partial
ordering on $P(f)$ defined as follows: for two elements $p$ and $q$ of
$P(f)$, $p < q$ if the boundary of the Fatou component containing $q$
separates $p$ from $x$.

Let $y$ be any minimal element with respect to this ordering.  Then
$y$ is a superattracting periodic point of period $p \geq 1$.  Let
$\Omega$ be the Fatou component containing $y$.  Let $E=P(f)-\{y\}$.
Since $y$ is minimal and $f$ is hyperbolic, $E \subset \Ins
(\gamma_{0})$, for some unique $\gamma_{0} \in \bdry \Omega$.
Moreover, $(f^{-p}E) \cap E \neq \vac$.  (If $p > 1 $ this is obvious,
since $f^{p}$ must fix every point in the orbit of $y$; if $p=1$, this
follows since not all points in the postcritical set can land on $y$
under one iterate of $f$.)  Hence $f^{-p}(E) \cap \Ins
(\gamma_{0}) \neq \vac$.  Since $P(f)=P(f^{p})$, the Jordan domain $\Out
(\gamma_{0}) $ contains a unique critical value of $f^{p}$ in its
closure, since $y$ is minimal and $f$ is hyperbolic.  Hence every
component of the preimage of $\Out (\gamma_{0})$ under $f^{p}$ is a
Jordan domain.  Proposition
\ref{prop:Omega fixed plus disc preimage}, {\it $\Omega$ fixed plus
disc preimage}, applied to $f^{p}$ now shows that $\Omega$ is a Jordan
domain.  Since $f$ is hyperbolic and $P(f^{p})=P(f)$, there are no
elements of $P(f)$ in $\partial \Omega $.  So every component $\Omega
'$ of $f^{-n}(\Omega), n \geq 0$, is also a Jordan domain, since
$\overline{\Omega '}$ is a branched cover of $\overline{\Omega }$
branched over at most one point which lies in the interior of $\Omega$.
\qed

{\bf Proof of Corollary 3 } By the above Theorem, the unique periodic
cycle of Fatou components of $f$ consists of Jordan domains.  Since
there are no critical points in the Julia set, there are no critical
values for iterates of $f$ in the boundaries of these Jordan domains.
Hence they all pull back to Jordan domains under iterates of $f$.
\qed

\subsection{Proof of Theorem 1}

If $f(z)$ is conjugate to $z^{d}$, then it is well-known that
$J(f)=S^{1}$.

Now suppose that upon conjugating by an automorphism of $\rs$, $f$ is
equal to a polynomial.  Let $\Omega$ be the basin of infinity.  Then
$J(f)=\bdry \Omega$.  If $\Omega$ is a Jordan domain, $J(f)$ is a
Jordan curve.  It then follows that there are exactly two Fatou
components $\Omega$ and $\Omega'$, and these components satisfy
$f^{-1}(\Omega)=\Omega$, $f^{-1}(\Omega')=\Omega'$.  Since $f$ is
critically finite, this implies that $f$ is of the form $z \mapsto
z^{d}, d \geq 2$.  If $\Omega$ is not a Jordan domain, then any other Fatou
component of $f$ is a component of $\rs - \cl{\Omega}$. By
Proposition \ref{prop:isjcurve}, these components are all Jordan domains.

So we may assume that $f$ is not conjugate to a polynomial or to a map
of the form $z \mapsto z^{d}$.  If there are no periodic critical
points, then $J(f) = \rs$, and there is nothing to prove.  Otherwise,
there is a periodic critical point $c_{1}$.  Since there are exactly
two critical points, each has multiplicity $d-1$, where $d=\deg(f)$.
If the period of $c_{1}$ is equal to one, then $f$ is conjugate to a
polynomial.  Hence we may assume $p \geq 2$.

Let $v_{1}$ be the image of $c_{1}$.  Let $\Omega_{0}$ be the Fatou
component containing $v_{1}$.  It suffices to show that $\Omega_{0}$
is a Jordan domain whose closure contains exactly one critical value
$v_{1}$.  For from this it follows that every Fatou component contains
at most one critical value in its closure.  Since every Fatou
component is eventually periodic, every Fatou component $\Omega'$ of
$f$ is a covering of a Jordan domain, branched over at most one point
in its closure, and hence $\Omega'$ is a Jordan domain.

Let $\Omega_{i}=f^{p-i}\Omega_{0}, i=1,...,p$; note that
$\Omega_{p}=\Omega_{0}$.  Since the $\Omega_{i}$ are Fatou components,
they are contained in unique components of $\rs -
\cl{\Omega_{0}}, i=1,...,p-1$.  By Proposition \ref{prop:fixed iff no
sign changes}, {\it $\Omega$ fixed iff no sign changes}, since
$\Omega$ is not fixed, there are sign changes in the set of sequences
$S_{out}(\Omega_{0})$.  Let $a_{\Omega_{0}}$ be the component of $\rs
- \cl{\Omega_{0}}$ containing $f(\Omega_{0})$, and let $\gamma_{0}$ be
the positively oriented boundary of $a_{\Omega_{0}}$.

\begin{figure}
\begin{center}~
\psfig{file=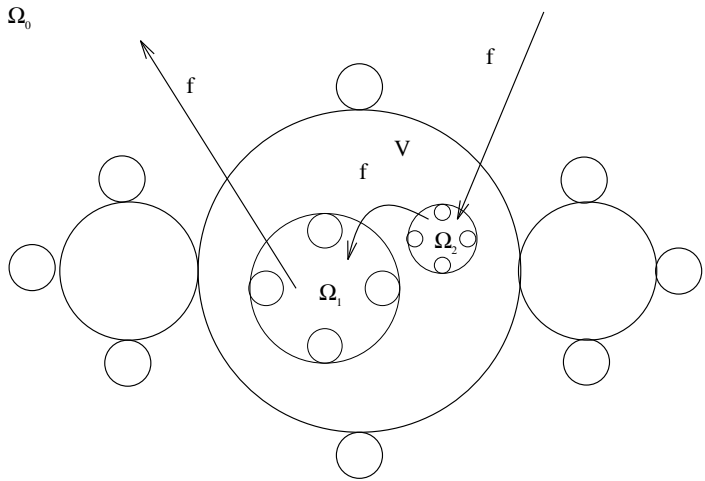,height=2.5in}
\end{center}
\caption{ The $\Omega_{i}, i=1,...,p-1$ are contained in
$V=a_{\Omega_{0}}$. Here $p=3$. \label{figure:basic_picture} }
\end{figure}

Since there are exactly two critical points, there are exactly two
critical values $v_{1}, v_{2}$.  Since $v_{1} \in \Omega_{0}$, a
Jordan domain in $\rs - \Omega$ contains at most one critical value
$v_{2}$ in its closure, hence the preimage under $f$ of every Jordan
domain in $\rs - \Omega_{0}$ is again a Jordan domain.  By Proposition
\ref{prop:Omega not fixed plus disc preimages}, {\it $\Omega$ not
fixed plus disc preimages}, we have that $\cup_{i=1}^{p-1}\Omega_{i}
\subset a_{\Omega_{0}}$.  We then have a basic picture of part of the
dynamics; see Figure
\ref{figure:basic_picture}.

Next, we prove
\gap

\noindent {\it If $v_{2} \in \Int (a_{\Omega_{0}})$, then
$\Omega_{0}$ is a Jordan domain}.
\gap

Let $D_{0}=\Out (\gamma_{0}) $.  Let $E=\{c_{1}\}$.  Then $E \subset
a_{\Omega_{0}}$ and $f^{p}(E)=E$. By Proposition \ref{prop:Omega
fixed plus disc preimage}, {\it $\Omega$ fixed plus disc preimage},
it suffices to prove that the unique component
$D_{p}$ of the preimage of $D_{0}$ under $f^{p}$ which contains
$\Omega_{0}$ is a Jordan domain.  We prove this by pulling back
$D_{0}$ along the orbit of $\Omega_{0}$ and using induction.

Let $D_{i}$ be the component of the preimage of $D_{0}$ under $f^{p}$
containing $\Omega_{i}, i=0,...,p$.  We first claim that $D_{1}
\subset a_{\Omega}$.  Since $v_{2} \in a_{\Omega_{0}}$, $D_{0}$
contains exactly critical value in its closure, so $D_{1}$ is a Jordan
domain.  By Part 1 of Proposition \ref{prop:Omega not fixed plus disc
preimages}, {\it $\Omega$ not fixed plus disc preimages}, we must have
$\gamma_{1} =\bdry D_{1} \subset a_{\Omega_{0}}$ and its sign must be negative.
Since the sign of $\gamma_{1}$ is negative, $\Out (\gamma_{1})=D_{1}
\subset a_{\Omega_{0}}$.

We now use induction.  Assume $D_{i}$ is a Jordan domain contained in
$\rs - \Omega_{0}, i=1,...,n<p$. Then $D_{n+1}$ is also a Jordan
domain since $D_{n}$ contains at most one critical value in its
closure.  A sign change between $\gamma_{n}$ and $\gamma_{n+1}$
implies that $\gamma_{n+1} \subset a_{\Omega_{0}}$, by Proposition
\ref{prop:Omega not fixed plus disc preimages}, {\it $\Omega$ not
fixed plus disc preimages}, and hence that $D_{n+1} \supset
D_{0}=f^{n+1}(D_{n+1})$.  But this implies by Proposition \ref{inF},
{\it Montel's theorem}, that $\Omega_{0}$ is a
Jordan domain fixed under the $(n+1)$st iterate of $f$, which is
impossible if $n+1 < p$.  The absence of a sign change then implies
that $D_{n+1} \subset \rs - \Omega_{0}$, and so the induction
proceeds.  Hence $D_{p-1}$ is a Jordan domain in $\rs -
\Omega_{0}$, and so $D_{p}$ is a Jordan domain.
\qed

So to prove Theorem \ref{thm:mainthm}, it suffices to prove (using the
notation in the preceding discussion)

\begin{prop}
\label{prop:v2}
The critical value $v_{2}$ is contained in the open disc
$a_{\Omega_{0}}$.
\end{prop}

As an immediate consequence, we have the following corollary:

\begin{cor}
Let $f$ be a critically finite rational map which has exactly two
critical points, and which is not conjugate to a polynomial.
Then no Fatou component of $f$ contains two critical values in its closure.
\end{cor}

\noindent {\bf Proof of Proposition:} We argue by contradiction.  Let
$D_{0}=\Out (\gamma_{0})$.  Let $D_{i}$ be
as above.  We will show that $\bdry D_{0} \subset \bdry D_{p}$, and
that $f^{p}|_{\bdry D_{0}}: \bdry D_{0} \rightarrow \bdry D_{0}$ is a
homeomorphism.  Since postcritically finite maps are expanding on
their Julia sets with respect to the canonical orbifold metric, by
Proposition \ref{prop:cf maps expand}, any compact connected set in
$J(f)$ mapped homeomorphically onto itself is a point.  This gives a
contradiction.

In order to carry out the argument, we need to show that $\bdry D_{i}
\subset \bdry \Omega_{i} \subset a_{\Omega_{0}}, i=0,1,...,p$.
This will be implied by the following lemma.  (We will only need the
case where $X_{0}$ is a Fatou component homeomorphic to an open disc
and $Y_{0}$ is homeomorphic either to the sphere minus a finite union of
disjoint closed discs, or to the sphere minus a finite union of closed
discs whose boundaries meet in exactly one common point to all of them.)

\begin{figure}
\begin{center}~ 
\psfig{file=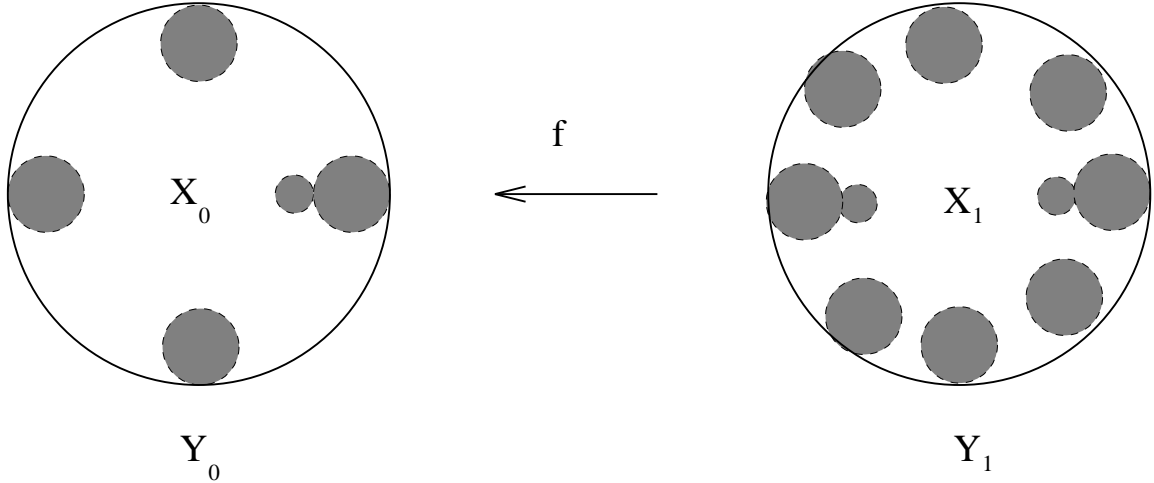,height=2.5in} 
\end{center}
\caption{$Y_{0}$ and $Y_{1}$ are the  large  discs.  $X_{i}$ is the complement
of the shaded discs in  $Y_{i}, i=1,2$.  \label{figure:boundaries} }
\end{figure}

\begin{lemma}
\label{boundaries}
Let $f: \rs \rightarrow \rs$ be a rational map.  Let $X_{0}$ and
$Y_{0}$ be proper open subsets of \rs   with $X_{0} \subset Y_{0}$.
Suppose $\partial Y_{0} \subset \partial X_{0}$.

\begin{enumerate}

\item If $Y_{1}=f^{-1}Y_{0}$ and $X_{1}=f^{-1}X_{0}$, then
$\partial Y_{1} \subset \partial X_{1}$. (See Figure
\ref{figure:boundaries}.)

\item If $Y_{1}$ is a component of $f^{-1}Y_{0}$, if $f |
_{\overline{Y_{1}}}: \overline{Y_{1}} \rightarrow \overline{Y_{0}}$ is
a homeomorphism, and if $X_{1}=(f |
_{\overline{Y_{1}}})^{-1}(X_{0})$, then $\partial Y_{1} \subset \partial
X_{1}$.

\end{enumerate}

\end{lemma}

{\bf Proof of Lemma:}

1.  Since $f$ is a nonconstant rational map, it is an open map, and so
for any proper open subset $Z \subset \rs $, $f^{-1}\partial Z =
\partial f^{-1}Z$.  So $\partial Y_{1}=\partial f^{-1}Y_{0} = f^{-1}
\partial Y_{0} \subset f^{-1} \partial X_{0} = \partial f^{-1}X_{0} =
\partial X_{1}$.

2.  Since $Y_{1}$ is a component of $f^{-1}Y_{0}$, $f(\partial Y_{1})
\subset \partial Y_{0}$.  Since $f: \overline{Y_{1}} \rightarrow
\overline{Y_{0}}$ is a homeomorphism and $X_{1}= (f |
_{\overline{Y_{1}}})
^{-1}(X_{0}), \partial X_{1} = (f | \overline{Y_{1}})^{-1}(\partial
X_{0})$.
Hence $\partial Y_{1} \subset \partial X_{1}$.
\qed

{\bf Remark:} The second statement is no longer true if we drop the
requirement that $\overline{Y_{1}}$ maps homeomorphically to
$\overline{Y_{0}}$.  For example, let $Y_{0}$ be the open unit disc,
let $X_{0}$ be the open disc minus the interval [0,1), and let
$f(z)=z^{2}$ map the Riemann sphere to itself.  Let $X_{1}$ be the
intersection of the upper half-plane $\{z| Im(z) > 0 \}$ with the unit
disc, and let $Y_{1}$ be the unit disc again.  Then all other
hypotheses of the lemma are satisfied but $\partial Y_{1} \not\subset
\partial X_{1}$.

{\bf Proof of Proposition \ref{prop:v2}, continued:}

Suppose $v_{2} \not \in \Ins (\gamma_{0})=\Int(a_{\Omega_{0}})$, and
let $D_{1}=f^{-1}(D_{0})$.

\begin{enumerate}

\item  We first show that $\partial D_{p} \subset \partial
\Omega_{0}$.

Since $v_{1} \in \Omega_{0} \subset D_{0}$, $D_{1}=f^{-1}(D_{0})$.  The
first case of Lemma \ref{boundaries} then applies, and so $\partial
D_{1} \subset \partial \Omega_{1}
\subset a_{\Omega_{0}}$.  It follows that $D_{1}$ must be contained in
$a_{\Omega_{0}}$.  For otherwise, $D_{1} \supset D_{0}$, and so $D_{1}
\subset F(f)$, by Proposition \ref{inF}, {\it Montel's theorem}.  But
then $p=1$ and so $\Omega$ is
fixed.

When $v_{2} \not\in a_{\Omega_{0}}$, the region $D_{1}$ is
homeomorphic to the complement of $d$ disjoint closed discs. If $v_{2}
\in \gamma_{0}=\bdry a_{\Omega_{0}}$, $D_{1}$ is homeomorphic to the
complement in the sphere of a union of $d$ closed discs whose
boundaries meet at exactly one point.  Note that in both cases the
boundary of $D_{1}$ consists of exactly $d$ lifts of $\gamma_{0}$,
each of which maps homeomorphically to $\gamma_{0}$.  This can be seen
as follows: $\rs - \cl{D_{0}}$ is a Jordan domain which contains
either no critical points in its closure, or one critical point in its
boundary.  Hence every component of its preimage is a Jordan domain,
there are exactly $d$ such components, and the boundary of each maps
homeomorphically under $f$.

We now argue by induction.  Assume for $1 \leq i < p$ that $D_{i}$ is
contained in $a_{\Omega_{0}}$, and that $\partial D_{i}
\subset \partial \Omega_{i}$.  We show that this
implies $\partial D_{i+1} \subset \partial \Omega_{i+1}$ if $i<p$, and
that $D_{i+1}$ is contained in $a_{\Omega_{0}}$ if $i<p-1$.  Since
$\Omega_{i} \subset a_{\Omega_{0}}$, and $D_{i} \subset
a_{\Omega_{0}}$ with $\partial D_{i} \subset \partial \Omega_{i}$,
$D_{i}$ is contained in $a_{\Omega_{0}}$, for otherwise $D_{i}$
contains $D_{0}$, implying by Proposition \ref{inF}, {\it Montel's
theorem},  that $\Omega_{0}$
is fixed under $f^{i+1}$.  The set $V=\Int(a_{\Omega_{0}})$ is a
Jordan domain containing at most one critical value in its closure.
Let $V'$ be the unique component of $f^{-1}V$ whose closure contains
$D_{i+1}$.  Since $V$ contains no critical values and $\bdry V$
contains at most one critical value, $f |
_{\overline{V'}}: \cl{V'} \rightarrow \cl{V}$ is a homeomorphism.  By
restriction, $f |_{\overline{D_{i+1}}}: \overline{D_{i+1}} \rightarrow
\overline{D_{i}}$ is also a homeomorphism.  We may now apply the
second case of Lemma \ref{boundaries} to conclude that $\partial
D_{i+1} \subset \partial \Omega_{i+1}$.  Moreover, if $i+1 < p$, then
$\Omega_{i+1} \subset a_{\Omega_{0}}$, and hence $\partial
\Omega_{i+1} \subset a_{\Omega_{0}}$.  $D_{i+1} \subset
a_{\Omega_{0}}$ if $i<p-1$, for otherwise $D_{i+1} \supset
D_{0}=f^{i+1}(D_{i+1})$.  Proposition \ref{inF}, {\it Montel's
theorem},  would then imply that
$D_{i+1} \subset F(f)$, whence $D_{i+1}=D_{0}=\Omega_{0}$ and
$\Omega_{0}$ is periodic of period strictly less than $p$, a
contradiction.

Hence $\bdry D_{p} \subset \bdry \Omega_{p}=\Omega_{0}$, and $D_{p}
\supset \Omega_{0}$.

\item We next claim that the boundary of every component of the
complement of $D_{i+1}$ (for convienience, let us call these {\it
boundary pieces} of $D_{i}$), $i=0,...,p-1$ maps injectively onto its
image under $f$, and so that every boundary piece of $D_{p}$ maps
homeomorphically onto its image $\partial D_{0}$ under $f^{p}$.
That this holds for $i=0$ has already been proved.  For
$1 \leq i \leq p-1$, the map $f: \overline{D_{i+1}} \rightarrow
\overline{D_{i}}$ is a homeomorphism, by the argument given in the
previous paragraph.

\item  We now claim that $f^{p}$ maps $\gamma_{0}$ homeomorphically
to itself, which we have shown is impossible.

By Proposition \ref{prop:Omega fixed}, {\it $\Omega$ fixed}, applied
to $f^{p}$, there must be some boundary piece of $D_{p}$ contained in
$a_{\Omega_{0}}$.  But since $\bdry D_{p} \subset \bdry \Omega$, this
implies that some boundary piece of $D_{p}$ is actually equal to
$\bdry a_{\Omega_{0}}=\gamma_{0}$.  The map $f$ restricted to a single boundary
component of $D_{i}, i=1,...,p$ is a homeomorphism, by the previous step,
and so the map $f^{p}$ sends $\gamma_{0}$ homeomorphically to itself.
\end{enumerate}
\qed

\section{Examples}
\label{section:examples}

\subsection{Examples where the Invariance Condition fails}

\begin{figure}
\begin{center}~
\psfig{file=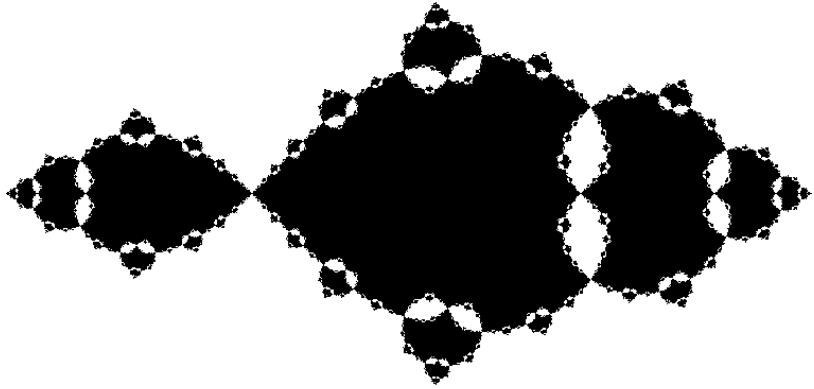,height=2in}
\end{center}
\caption{The degree 4 pseudo-basilica \label{pseudobasilica4} }
\end{figure}

J. Kahn, C. McMullen and the author discovered a degree four map
where the Invariance Condition appears to fail: the map turned out to
be
\[  z \mapsto \frac{3(z-1)^3(z+3)}{3-8z + 6z^{2}} 	\]
whose Julia set is given in Figure \ref{pseudobasilica4}.  The black
regions just to the right and left of the pinched point are the
immediate basins of attraction of a period two superattracting cycle.
All black regions eventually map onto these basins.  The white regions
all eventually map onto the basin of infinity, which is forward-invariant.

This map is one member in a family of maps of varying degree:  set
\[	f_{d}(z) =  N_{d} \circ p_{d} \circ M(z)	\]
where $M(z)=\frac{z-1}{z}$, $p_{d}(z)=(d-1)z^{d}-dz^{d-1} +1$, and
$N_{d}=(1-d)\frac{z-1}{z}$.  For this family, infinity is a simple
critical point, 1 is a critical point of local degree $d-1$, 0 is a
critical point of local degree $d$, 1 maps to 0, 0 maps to $1-d$, and
$1-d$ maps to 0.  For $d=2$ one obtains a map conjugate to $z \mapsto
z^{2}-1$, whose Julia set is called the basilica.  For $d \geq 3$,
however, on the basin of the unbounded Fatou component $\Omega$, the
map is conjugate to $z \mapsto z^{2}$, but $\bdry \Omega$ appears to
homeomorphic to a figure-8.  We refer to the Julia sets of $f_{d}$ as
``pseudo-basilicas''.

We now give a direct argument in the case $d=3$ which shows that the
Invariance Condition fails for the boundary of the basin of infinity.
For the definition of Thurston equivalence of branched
coverings, and Thurston's theorem on the existence of a rational map
in a given Thurston class, see e.g. \cite{rees:parami} or \cite{DH1}.

\subsection{The degree three pseudobasilica}

\begin{figure}
\begin{center}~
\psfig{file=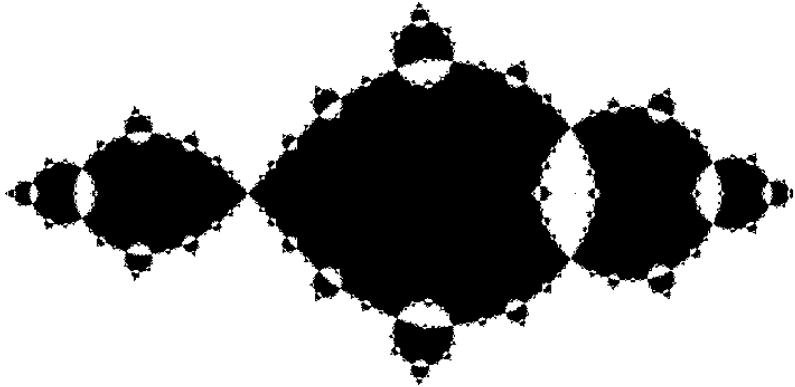,height=2in}
\end{center}
\caption{The pseudo-basilica in degree 3\label{pseudobasilica3} }
\end{figure}

Let $g$ be the map $f_{3}(z)=\frac{(z+2)(z-1)^{2}}{\frac{3}{2}z-1}$
and $ \Omega $ be the basin of infinity of $g$.  Figure
\ref{pseudobasilica3} is a picture of its Julia set.  The point at
which $\bdry \Omega$ appears pinched is not a critical point of $g$.
We give a proof that $g|_{\bdry \Omega}$ fails the invariance
condition which depends strongly on the degree and the fact that the
map is real.

Since the point at infinity is a simple superattracting critical point
and $g$ is postcritically finite, a theorem of B\"{o}ttcher
(\cite{MIL1}, Theorem 6.7) implies that there exists a unique Riemann
map $\phi: (\Delta, 0) \rightarrow (\Omega,
\infty)$ such that $\phi(z^{2}) = f(\phi(z))$.  In what follows,
``ray'' means a ray for the map $\phi$.

{\bf Step 1} We first claim that the 1/3 and 2/3 rays land at a common
fixed point.   First, we show that every point of period less than or equal
to two is real.  There are ten total, one of which is the point at
infinity. The remainder are the finite solutions of the equation
$g^{2}(z)-z=0$. Two solutions are the period 2 attractors 0 and -2.
Another is the fixed point 2, which is the landing point of the 0 ray
in the basin of infinity. The remaining six solutions are roots of the
polynomial $(4z^{4}-2z^{3}-15z^{2}+16z-4)(2z^{2}+z-2)$, all of which
are real.  We next claim that the 1/3 and 2/3 rays in $ \Omega $
land at a common fixed point $p$.  Since $g$ is real, it commutes with
conjugation, so $R_{1/3}=\overline{R_{2/3}}$.  The landing point of
$R_{1/3}$ is therefore the complex conjugate of the landing point of
$R_{2/3}$.  But these landing points are points of period less than or
equal to two, so by Step 1, they must be real and hence equal.  Since
the two rays are exchanged under the dynamics, the common landing point
$x$ is actually fixed under $g$.

{\bf Step 2} Let $C$ be the closed curve which is the union of the
point at infinity, the 1/3 and 2/3 rays, and $p$.  We claim that $C$
separates 0 from -2 in \rs .  For otherwise, one component of $U$ of
$\rs - C$ is an open disc containing no elements of $P(g)=P(g^{2})$.
The preimages of $U$ under $g^{2}$ are then all disjoint open discs.
The curve $C$ is fixed as an oriented curve under $g^{2}$.  So for
some preimage $V$ of $U$ under $g^{2}$, $V=U$, and so $U$ must be
contained in the Fatou set.  But this then implies that $C$ cannot
separate $\partial \Omega $, contradicting Theorem
\ref{rays_separate}.

{\bf Step 3} We next claim that there is a unique (up to combinatorial
equivalence) real degree three branched cover of the sphere with the
same postcritical data.  Since any critically finite branched covering
$G$ with $|P(G)|=3$ and hyperbolic orbifold is Thurston equivalent to
a rational map, it suffices to prove there is a unique rational map
with this data.  By conjugating we may assume that infinity is a fixed
simple critical point, 1 maps with local degree two to zero and -2
maps with local degree one to zero.  These conditions imply that $g$
is of the form $g(z)= \frac{(z-1)^{2}(z+2)}{az-b}$.  If we require
that zero is to map with local degree three onto its image and then
back to itself, there are a unique parameters $a$ and $b$, namely
$a=3/2$ and $b=-1$.

\begin{figure}
\begin{center}~
\psfig{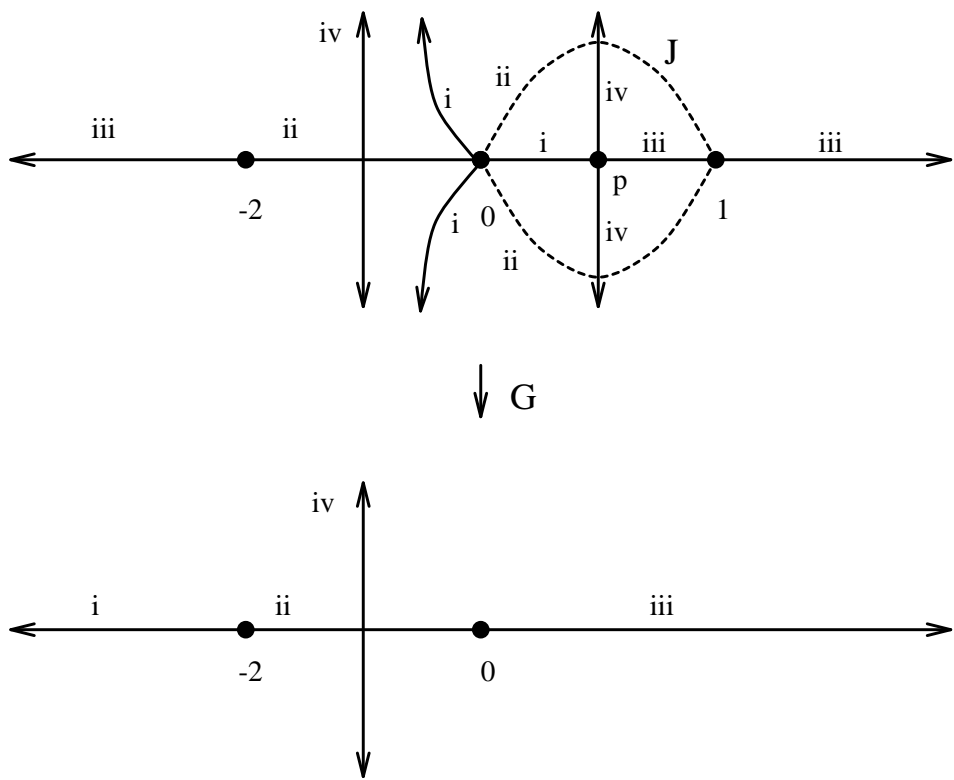}
\end{center}
\caption{The branched covering $G$.  The bounded region formed by the
dashed curve $J$ maps onto the complement of the arc numbered by
(ii). \label{fig:semibasex}}
\end{figure}

{\bf Step 4}  The previous step implies that the branched
cover $G$ described in Figure \ref{fig:semibasex}
is Thurston equivalent to the map $g$.  The top figure is to be
overlaid the bottom one to form a critically finite branched covering
of the sphere to itself.

For this map, a loop $\gamma$, separating 0 from -2,
and represented by the line (iv) in the bottom half of Figure
\ref{fig:semibasex} union the point at infinity, has two preimages.
One preimage is homotopic to $\gamma$ and maps by degree $-1$, while
the other is a closed curve mapping to $\gamma$ by degree two,
represented as the union of the point at infinity together with the
two arcs in the top figure labelled (iv) passing through the pole $p$.

Since $G$ is Thurston equivalent to $g$, and $C$ corresponds to the
curve $\gamma$ under the obvious Thurston equivalence, it follows that the
curve $C$ must also have a preimage which maps to $C$ by degree two.
Hence the endpoints of $R_{1/6}$ and $R_{5/6}$ are necessarily
distinct, implying that $p$ has three preimages in $\partial \Omega $
under $g$.  For any polynomial $p(z)$ of degree $d$ with basin of
infinity $\Omega'$, $J(p)=\bdry \Omega'$, and so $\bdry \Omega'$ is
totally invariant under $p$.  If $g|_{\bdry \Omega}$ were topologically
conjugate to a polynomial, a generic point in $\bdry \Omega$ would
then have three preimages, and so $\bdry \Omega$ would be totally
invariant under $g$.  Since $g(\Omega)=\Omega$, if $\bdry \Omega$ is
totally invariant, then so is $\Omega$, which it is not.  Hence
$g|_{\bdry \Omega}$ cannot be topologically conjugate to the dynamics
of any polynomial on its Julia set.

\rmk A similar proof works to show that the
maps $f_{d}$ possess the same property.  The only significant
difference is that a different argument in Step 1 is required.  One
can prove this using the fact that $\gamma$ is fixed up to homotopy
relative to $P(G)$ as an unoriented curve under $G$, together with the
fact that the maps $f_{d}$ are expanding on their Julia sets.

\rmk For the map $g$, if $\Omega$ is the basin of infinity, the set
$A$ consists of precisely two elements.  Let $a_{-2}, a_{0}$ be the closures
of the components of $\rs - \cl{\Omega}$ containing $-2$ and $0$
respectively.  Let $\gamma_{-2}=\bdry a_{-2}$ and $\gamma_{0}=\bdry
a_{0}$ with positive orientation relative to $\Omega$.  Then
$\gamma_{-2}$ has a unique lift $\widetilde{\gamma_{-2}}$ mapping by
degree $+3$.  There is an open topological arc $\alpha \subset
\widetilde{\gamma_{-2}}$ mapping homeomorphically to
$\gamma_{-2}-\{x\}$; this open arc is the portion of the boundary of
the immediate basin of $0$ lying between the landing points of the
$1/6$ and $5/6$ rays; see Figure \ref{pseudobasilica3}.  Thus
$\widetilde{\gamma_{-2}} \subset a_{0}$ is not contained in $\bdry \Omega$.
The curve $\gamma_{0}$ has two lifts. One is $\gamma_{-2}$ which maps
by degree $+1$; the other maps by degree $+2$ and is contained in
$a_{0}$ but not in $\bdry a_{0}$.

\subsection{Other interesting examples}

\begin{figure}

\begin{center}~
\psfig{file=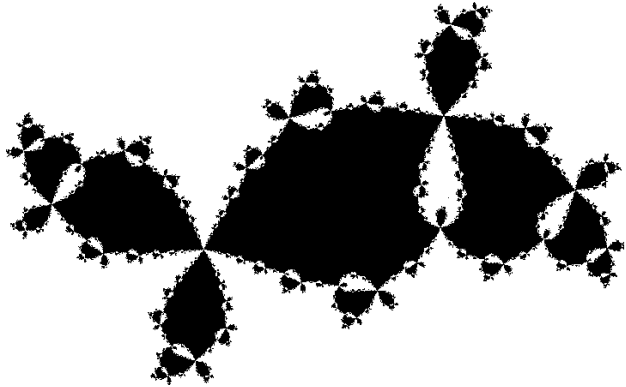,height=2in}
\end{center}
\caption{A degree 3 pseudo-rabbit \label{pseudorabbit3} }

\end{figure}

The family $g_{r}(z)=\frac{r}{d-1} N_{d} \circ p_{d} \circ M $,
provide other examples of maps which appear to fail the Invariance
Condition.  For example, if $d=3$, and if $r$ is a
solution of $g_{r}^{3}(0)=0, g_{r}(0) \neq 0$, the Julia sets look
quite interesting.  A complex solution is $r \approx 1.34781+1.02885i$
and yields a ``pseudorabbit'', as is shown in Figure
\ref{pseudorabbit3}.

\end{document}